\documentclass[preprint,12pt]{elsarticle}

\usepackage{listings}
\usepackage{color}

\definecolor{dkgreen}{rgb}{0,0.6,0}
\definecolor{gray}{rgb}{0.5,0.5,0.5}
\definecolor{mauve}{rgb}{0.58,0,0.82}

\usepackage{tikz}
\usetikzlibrary{backgrounds}
\usepackage{graphicx}
\usepackage{bm}
\usepackage{amsmath}
\usepackage{setspace}
\usepackage{hyperref}
\usepackage{subfigure}
\usepackage{multirow}
\hypersetup{
    colorlinks=true,       % false: boxed links; true: colored links
}
\usepackage{amssymb}
\usepackage{lineno}
\usepackage{diagbox}

\journal{Journal Name}

\begin{document}

\begin{frontmatter}


\title{Multistep Neural Networks for Data-driven Discovery of Nonlinear Dynamical Systems}



\author{Maziar Raissi$^{1}$, Paris Perdikaris$^{2}$, and George Em Karniadakis$^{1}$}
\address{$^{1}$Division of Applied Mathematics, Brown University,\\ Providence, RI, 02912, USA\\
$^{2}$Department of Mechanical Engineering and Applied Mechanics,\\ University of Pennsylvania,\\ Philadelphia, PA, 19104, USA}
%\address{Division of Applied Mathematics, Brown University,\\ Providence, RI, 02912, USA}

\begin{abstract}

% blend classical ideas from time-stepping with nonlinear regression tools for systems identification

The process of transforming observed data into predictive mathematical models of the physical world has always been paramount in science and engineering. Although data is currently being collected at an ever-increasing pace, devising meaningful models out of such observations in an automated fashion still remains an open problem. In this work, we put forth a machine learning approach for identifying nonlinear dynamical systems from data. Specifically, we blend classical tools from numerical analysis, namely the multi-step time-stepping schemes, with powerful nonlinear function approximators, namely deep neural networks, to distill the mechanisms that govern the evolution of a given data-set. We test the effectiveness of our approach for several benchmark problems involving the identification of complex, nonlinear and chaotic dynamics, and we demonstrate how this allows us to accurately learn the dynamics, forecast future states, and identify basins of attraction. In particular, we study the Lorenz system, the fluid flow behind a cylinder, the Hopf bifurcation, and the Glycoltic oscillator model as an example of complicated nonlinear dynamics typical of biological systems.

\end{abstract}

\begin{keyword}

Machine learning \sep Systems identification \sep Reduced order modeling \sep Data-driven forecasting \sep Lorenz system \sep Navier-Stokes


\end{keyword}

\end{frontmatter}

% \linenumbers

\section{Introduction}

% dynamical systems
Dynamical systems play a key role in shaping our understanding of the physical world and in defining our ability to predict the evolution of a given process. From the simple swinging motion of a clock pendulum to the complex flow around an airplane wing, the mathematical modeling of dynamical systems can yield a set of tools with which we can analyze the way the current state of the system depends on the past, and predict the possible states we may encounter in the future. Often such tools are precisely known, usually coming in the form of differential equations that are derived from first physical principles, such as the conservation of energy, mass, and momentum \cite{courant2008methods}. However, in many cases, the sheer complexity of a system can prohibit our complete understanding and render a first principles approach infeasible. In this setting, one may be only able to postulate crude and potentially overly simplified models based on a given a set of empirical observations (see e.g., models for tumor growth \cite{michor2004dynamics}, social dynamics \cite{castellano2009statistical}, and the stock market \cite{gavin1989stock}). In the present era of abundant data and advanced machine learning capabilities, a natural question arises: can we automatically discover sufficiently sophisticated and accurate mathematical models of complex dynamical systems directly from data?\\

% systems identification and discovery of governing equations
The answer to this question pertains to the well-established field of systems identification \cite{ljung1998system}. Discriminating between {\em white-, gray-, and black-box} approaches, depending on whether a first principles modeling approach is fully, partially, or not admissible, systems identification aims to devise mathematical models for predicting a future state of a system, given the evolution of a set of previously observed or latent states. Specifically, in the context of identifying nonlinear dynamics, there exist several deterministic and probabilistic tools including radial basis functions \cite{chen1990non}, neural networks \cite{cochocki1993neural}, Gaussian processes \cite{kocijan2005dynamic, raissi2017machine,raissi2017hidden,raissi2017inferring,raissi2017numerical}, and nonlinear auto-regressive models such as NARMAX \cite{billings2013nonlinear} and recurrent neural networks \cite{goodfellow2016deep}. A common theme among all such methods is the pursuit of learning a nonlinear and potentially multi-variate mapping $\bm{f}$ that predicts the future system states given a set of data describing the present and past states. More recently, approaches based on symbolic regression \cite{schmidt2009distilling}, sparse regression, and compressive sensing \cite{brunton2016discovering, rudy2017data} were able to go beyond estimating a {\em black-box} approximation of the dynamics given by $\bm{f}$, and return more interpretable models that can uncover the full parametric form of an underlying governing equation. However, in order to obtain sparse representation of the dynamics, the aforementioned approaches have to rely upon the nontrivial task of choosing ``appropriate" sets of basis functions. Consequently, investigating ways of incorporating broader function search spaces is an important area of current and future research.\\

% this work
In this work, we introduce a novel approach to nonlinear systems identification that combines the classical multistep family of time-stepping schemes from numerical analysis \cite{iserles2009first} with deep neural networks. Inspired by recent developments in {\em physics-informed deep learning} \cite{raissi2017physics_I,raissi2017physics_II}, we construct structured nonlinear regression models that can discover the  dynamic dependencies in a given set of temporal data-snapshots, and return a closed form model that can be subsequently used to forecast future states or identify basins of attraction. In contrast to recent approaches to systems identification \cite{brunton2016discovering, rudy2017data}, here we do not have to have direct access or approximations to temporal gradients because the time derivatives are discretized using classical time-stepping rules. Moreover, we are using a richer family of function approximators and consequently we do not have to commit to a particular class of basis functions such as polynomials or sines and cosines. This comes at the cost of losing interpretability of the learned dynamics. However, there is nothing hindering the use of a
particular class of basis functions and obtain more interpretable equations. \\

This paper is structured as follows. In section \ref{sec:setup} we provide a detailed overview of the proposed methodology.
In section \ref{sec:Cubic2D}, we investigate the performance of the proposed framework by applying our algorithm to the two-dimensional damped harmonic oscillator. We then explore the identification of chaotic dynamics of the Lorenz system in section \ref{sec:Lorenz}. As an example of a high dimensional dynamical systems, in section \ref{sec:NavierStokes}, we study the Navier-Stokes equations describing the fluid flow behind a cylinder. To illustrate the ability of our method to identify parameterized dynamics, we consider the Hopf normal form in section \ref{sec:Hopf}. As an example of complicated nonlinear dynamics typical of biological systems, we explore the glycolytic oscillator model in section \ref{sec:glycolytic}. It should be highlighted that all of the examples considered in this work are inspired by the pioneering work of Brunton \emph{et. al.} \cite{brunton2016discovering}. Moreover, all data and codes used in this manuscript are publicly available on GitHub at \url{https://github.com/maziarraissi/MultistepNNs}.

\section{Problem setup and solution methodology}\label{sec:setup}
Let us consider nonlinear dynamical systems of the form\footnote{It is straightforward to generalize the dynamics to include
parameterization, time dependence, and forcing. In particular, parameterization, time dependence, and external forcing or feedback control $\bm{u}(t)$ may be added to the vector field according to
\[
\dot{\bm{x}} = \bm{f}\left(\bm{x}, \bm{u}, t; \lambda\right),\ \ \ \dot{t} = 1,\ \ \ \dot{\lambda} = 0.
\]
}
\begin{equation}\label{eq:DynamicalSystems}
\frac{d}{d t} \bm{x}(t) = \bm{f}\left(\bm{x}(t)\right),
\end{equation}
where the vector $\bm{x}(t) \in \mathbb{R}^D$ denotes the state of the system at time $t$ and the function $\bm{f}$ describes the evolution of the system. Given noisy measurements of the state $\bm{x}(t)$ of the system at several time instances $t_1, t_2, \ldots, t_N$, our goal is to determine the function $\bm{f}$ and consequently discover the underlying dynamical system \eqref{eq:DynamicalSystems} from data. We proceed by applying the general form of a linear multistep method with  $M$ steps to equation \eqref{eq:DynamicalSystems} and obtain
\begin{equation}\label{eq:multistep}
\sum_{m=0}^M \left[\alpha_m \bm{x}_{n-m} + \Delta t \beta_m \bm{f}(\bm{x}_{n-m})\right] = 0, \ \ \ n = M, \ldots, N.
\end{equation}
Here, $\bm{x}_{n-m}$ denotes the state of the system $\bm{x}(t_{n-m})$ at time $t_{n-m}$. Different choices for the parameters $\alpha_m$ and $\beta_m$ result in specific schemes. For instance, the trapezoidal rule 
\begin{equation}\label{eq:trapezoidal}
\bm{x}_n = \bm{x}_{n-1} + \frac{1}{2} \Delta{t} \left(\bm{f}(\bm{x}_n) + \bm{f}(\bm{x}_{n-1})\right),\ \ \ n = 1, \ldots, N,
\end{equation}
corresponds to the case where $M = 1$, $\alpha_0 = -1$, $\alpha_1 = 1$, and $\beta_0 = \beta_1 = 0.5$. We proceed by placing a neural network prior on the function $\bm{f}$. The parameters of this neural network can be learned by minimizing the mean squared error loss function
\begin{equation}\label{eq:MSE}
MSE := \frac{1}{N-M+1}\sum_{n=M}^{N} |\bm{y}_n|^2,
\end{equation}
where
\begin{equation}\label{eq:multistep_NN}
\bm{y}_n := \sum_{m=0}^M \left[\alpha_m \bm{x}_{n-m} + \Delta t \beta_m \bm{f}(\bm{x}_{n-m})\right], \ \ \ n = M, \ldots, N,
\end{equation}
is obtained from the multistep scheme \eqref{eq:multistep}.

\section{Results}

\subsection{Two-dimensional damped oscillator}\label{sec:Cubic2D}
\begin{figure}[!t]
\includegraphics[width = 1.0\textwidth]{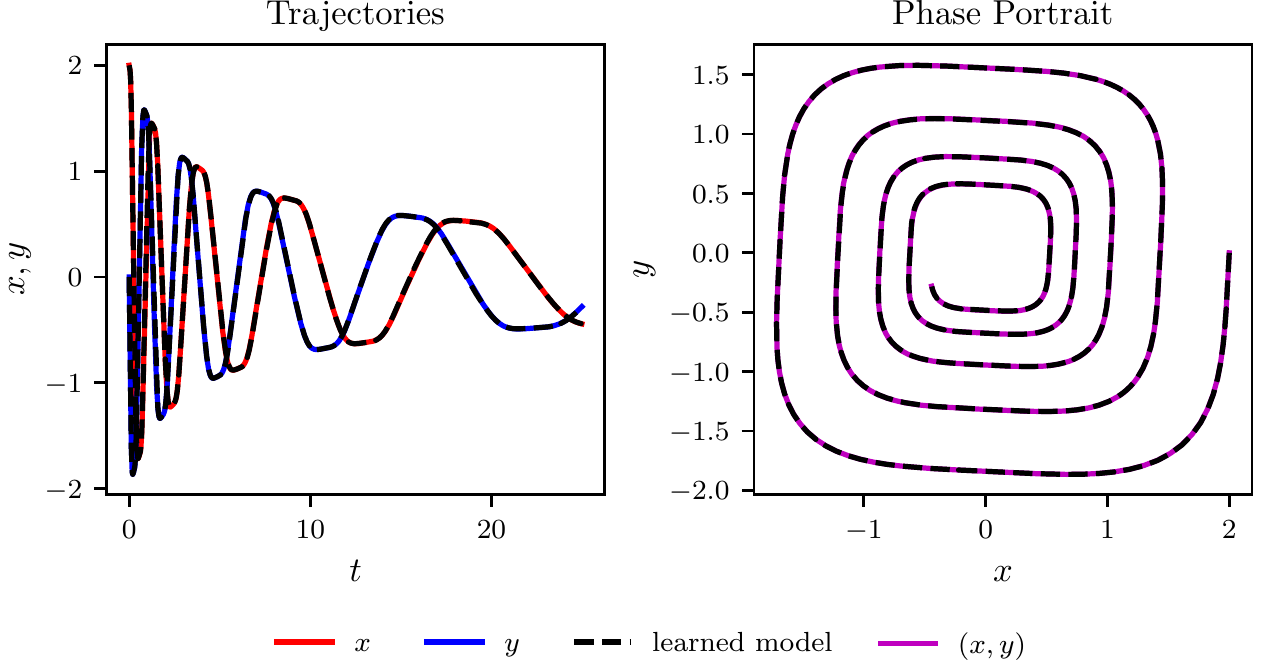}
\caption{{\em Harmonic Oscillator:} Trajectories of the two-dimensional damped harmonic oscillator with cubic dynamics are depicted in the left panel while the corresponding phase portrait is plotted in the right panel. Solid colored lines represent the exact dynamics while the dashed black lines demonstrate the learned dynamics. The identified system correctly captures the form of the dynamics and accurately reproduces the phase portrait.}
\label{fig:Cubic2D}
\end{figure}

As a first illustrative example, let us consider the two-dimensional damped harmonic oscillator with cubic dynamics; i.e.,
\begin{equation}
\begin{array}{l}
\dot{x} = -0.1\ x^3 + 2.0\ y^3,\\
\dot{y} = -2.0\ x^3 - 0.1\ y^3.
\end{array}
\end{equation}
We use $[x_0\ y_0]^T = [2\ 0]^T$ as initial condition and collect data from $t = 0$ to $t = 25$ with a time-step size of $\Delta t = 0.01$. The data are plotted in figure \ref{fig:Cubic2D}. We employ a neural network with one hidden layer and $256$ neurons to represent the nonlinear dynamics. As for the multistep scheme \eqref{eq:multistep} we use Adams-Moulton with $M=1$ steps (i.e., the trapezoidal rule). Upon training the neural network, we solve the identified system using the same initial condition as the one above. Figure \ref{fig:Cubic2D} provides a qualitative assessment of the accuracy in identifying the correct nonlinear dynamics. Specifically, by comparing the exact and predicted trajectories of the system, as well as the resulting phase portraits, we observe that the algorithm can correctly capture the dynamic evolution of the system. \\

To investigate the performance of the proposed work-flow with respect to different linear multi-step methods, we have considered the three  families that are most commonly used in practice: Adams-Bashforth (AB) methods, Adams-Moulton (AM) methods, and the backward differentiation formulas (BDFs). In tables \ref{tab:scheme_vs_M_x} and \ref{tab:scheme_vs_M_y}, we report the relative $\mathcal{L}_2$ error between trajectories of the exact and the identified systems for different members of the class of linear multi-step methods. 
Interestingly, the Adams-Moulton scheme seems to consistently return more accurate results compared to the Adams-Bashforth and BDF approaches. One intuitive explanation for this behavior stems from a closer inspection of equation \ref{eq:multistep}. Specifically, the arrangement of the resulting terms for the Adams-Moulton schemes leads to a higher throughput of training data flowing through the 
neural network during model training as compared to the Adams-Bashforth and BDF cases. This
helps regularize the neural network and eventually achieve a better calibration during training. Also, out
of the Adams-Moulton family, the trapezoidal rule seems to work the best in practice perhaps due to its superior stability properties \cite{iserles2009first}. These performance characteristics should be interpreted as product of empirical evidence, and not as concrete theoretical properties of the method. Identification of the latter requires more extensive systematic studies that go beyond the scope of this paper. \\

\begin{table}[!t]
\centering
\begin{tabular}{|l||ccccc|} 
\hline
\diagbox{Scheme}{M} & 1 & 2 & 3 & 4 & 5 \\ \hline\hline
Adams-Bashforth & 1.5e+00 & 3.1e-02 & 1.2e-01 & 4.3e-02 & 1.2e-02 \\
Adams-Moulton & 8.8e-03 & 1.2e-02 & 1.6e-02 & 6.3e-03 & 1.1e-02 \\
BDF & 1.3e+00 & 8.8e-03 & 1.3e-02 & 1.4e-02 & 1.7e-02 \\
\hline
\end{tabular}
\caption{{\em  Harmonic Oscillator:} Relative $\mathcal{L}_{2}$ error between the predicted and the exact trajectory for the first dynamic component $x(t)$ integrated up to time $t=25$, for different member families of the class of multistep methods, and different number of steps $M$. Here, the training data is assumed to be noise free, and the neural network architecture is kept fixed to have one hidden layer and 256 neurons.} \label{tab:scheme_vs_M_x}
\end{table}
\begin{table}[!t]
\centering
\begin{tabular}{|l||ccccc|} 
\hline
\diagbox{Scheme}{M} & 1 & 2 & 3 & 4 & 5 \\ \hline\hline
Adams-Bashforth & 1.5e+00 & 3.0e-02 & 9.7e-02 & 3.5e-02 & 1.2e-02 \\
Adams-Moulton & 8.8e-03 & 1.0e-02 & 1.6e-02 & 5.8e-03 & 1.1e-02 \\
BDF & 1.3e+00 & 8.6e-03 & 9.9e-03 & 1.4e-02 & 1.5e-02 \\
\hline
\end{tabular}
\caption{{\em  Harmonic Oscillator:} Relative $\mathcal{L}_{2}$ error between the predicted and the exact trajectory for the second dynamic component $y(t)$ integrated up to time $t=25$, for different member families of the class of multistep methods, and different number of steps $M$. Here, the training data is assumed to be noise free, and the neural network architecture is kept fixed to have one hidden layer and 256 neurons.} \label{tab:scheme_vs_M_y}
\end{table}

In tables \ref{tab:dt_vs_noise_x} and \ref{tab:dt_vs_noise_y}, we study the robustness of our results with respect to the gap $\Delta t$ between pairs of data and with respect to noise in the observations of the system. These results fail to reveal a consistent pattern as larger time-step sizes $\Delta t$ and larger noise corruption levels sometimes lead to superior accuracy and other times to inferior. In the latter cases, the reasons are obvious, namely that a larger gap $\Delta t$ makes the approximation in time less accurate, while too much noise is devastating because it would be harder to distinguish between noise and the true dynamics. On the other hand, we do observe some cases in which larger $\Delta t$ and noise levels may actually help. In these cases, we believe that input noise can act as a regularization mechanism that increases the robustness of the model training procedure, similar to how it has been previously proposed in the neural network literature (see for e.g., denoising autoencoders \cite{vincent2008extracting}). Along the same lines, a bigger temporal gap $\Delta t$ helps because it makes two consecutive time snapshots carry more information simply because they are more dissimilar to one another. On the contrary, if $\Delta t$ is too small, the importance of the neural network becomes less and less pronounced as seen in equation \eqref{eq:multistep}, hence model training becomes infeasible. These empirical results indicate that there exists a problem-dependent sweet spot for the admissible values of the time-step and noise levels that can lead to the best predictive accuracy. Although one may have no control over the noise corrupting the data, the temporal gap $\Delta t$ could be treated as another hyper-parameter like the number of neurons and hidden layers when setting up the neural network.\\

\begin{table}[!t]
\centering
\begin{tabular}{|l||ccc|} 
\hline
\diagbox{$\Delta{t}$}{noise} & 0.00\% & 0.01\% & 0.02\% \\ \hline\hline
0.01 & 5.7e-03 & 2.4e-02 & 2.2e-01 \\
0.02 & 1.8e-02 & 1.1e-01 & 1.3e-01 \\
0.03 & 3.8e-02 & 9.2e-02 & 7.8e-01 \\
0.04 & 5.4e-02 & 4.0e-02 & 8.7e-01 \\
0.05 & 8.3e-02 & 2.9e-01 & 9.2e-02 \\
\hline
\end{tabular}
\caption{{\em Harmonic Oscillator:} Relative $\mathcal{L}_{2}$ error between the predicted and the exact trajectory for the first dynamic component $x(t)$ integrated up to time $t=25$, for different noise magnitudes, and different gap $\Delta{t}$ between pairs of snapshots. Here, we are employing the trapezoidal time-stepping scheme and the neural network architecture is kept fix to have one hidden layer and 256 neurons.} \label{tab:dt_vs_noise_x}
\end{table}
\begin{table}[!t]
\centering
\begin{tabular}{|l||ccc|} 
\hline
\diagbox{$\Delta{t}$}{noise} & 0.00\% & 0.01\% & 0.02\% \\ \hline\hline
0.01 & 5.9e-03 & 2.1e-02 & 2.2e-01 \\
0.02 & 1.7e-02 & 1.0e-01 & 1.1e-01 \\
0.03 & 3.3e-02 & 9.0e-02 & 7.9e-01 \\
0.04 & 5.2e-02 & 3.8e-02 & 8.7e-01 \\
0.05 & 7.2e-02 & 2.6e-01 & 8.2e-02 \\
\hline
\end{tabular}
\caption{{\em Harmonic Oscillator:} Relative $\mathcal{L}_{2}$ error between the predicted and the exact trajectory for the first dynamic component $y(t)$ integrated up to time $t=25$, for different noise magnitudes, and different gap $\Delta{t}$ between pairs of snapshots. Here, we are employing the trapezoidal time-stepping scheme and the neural network architecture is kept fix to have one hidden layer and 256 neurons.} \label{tab:dt_vs_noise_y}
\end{table}

Finally, tables \ref{tab:layers_vs_neurons_x} and \ref{tab:layers_vs_neurons_y} study the robustness of our results with respect to the neural network structure. For this case, more accurate results seem to be obtained with increasing network depth, although increasing the network width seems to have a negative affect for more than $128$ neurons per layer. To fully quantify sensitivity with respect to network architecture a more systematic study involving multiple data-sets is needed.

\begin{table}[!t]
\centering
\begin{tabular}{|l||ccc|} 
\hline
\diagbox{Layers}{Neurons} & 64 & 128 & 256 \\ \hline\hline
1 & 9.8e-03 & 6.1e-03 & 3.7e-02 \\
2 & 3.6e-03 & 1.2e-02 & 2.4e-02 \\
3 & 3.4e-03 & 1.6e-02 & 4.2e-02 \\
\hline
\end{tabular}
\caption{{\em  Harmonic Oscillator:} Relative $\mathcal{L}_{2}$ error between the predicted and the exact trajectory for the first dynamic component $x(t)$ integrated up to time $t=25$, for different neural network architectures. Here, the training data is assumed to be noise free, the time step size is kept fixed at $\Delta t=0.01$, and the number of Adams-Moulton steps is fixed at $M=1$.}\label{tab:layers_vs_neurons_x}
\end{table}
\begin{table}[!t]
\centering
\begin{tabular}{|l||ccc|} 
\hline
\diagbox{Layers}{Neurons} & 64 & 128 & 256 \\ \hline\hline
1 & 7.2e-03 & 5.4e-03 & 3.5e-02 \\
2 & 3.3e-03 & 9.1e-03 & 2.0e-02 \\
3 & 3.0e-03 & 1.4e-02 & 3.7e-02 \\
\hline
\end{tabular}
\caption{{\em  Harmonic Oscillator:} Relative $\mathcal{L}_{2}$ error between the predicted and the exact trajectory for the second dynamic component $y(t)$ integrated up to time $t=25$, for different neural network architectures. Here, the training data is assumed to be noise free, the time step size is kept fixed at $\Delta t=0.01$, and the number of Adams-Moulton steps is fixed at $M=1$.}\label{tab:layers_vs_neurons_y}
\end{table}

\subsection{Lorenz system}\label{sec:Lorenz}
\begin{figure}[!t]
\includegraphics[width = 1.0\textwidth]{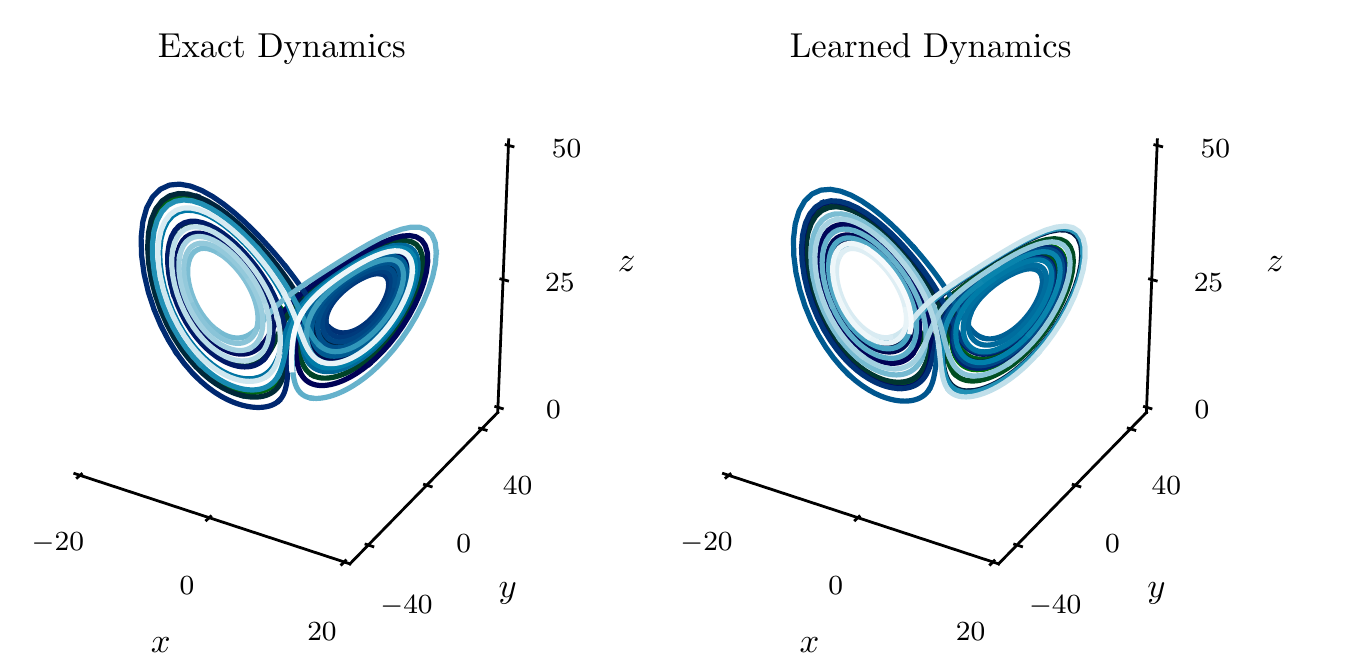}
\caption{{\em Lorenz System:} The exact phase portrait of the Lorenz system (left panel) is compared to the corresponding phase portrait of the learned dynamics (right panel).}
\label{fig:Lorenz}
\end{figure}

\begin{figure}[!t]
\includegraphics[width = 1.0\textwidth]{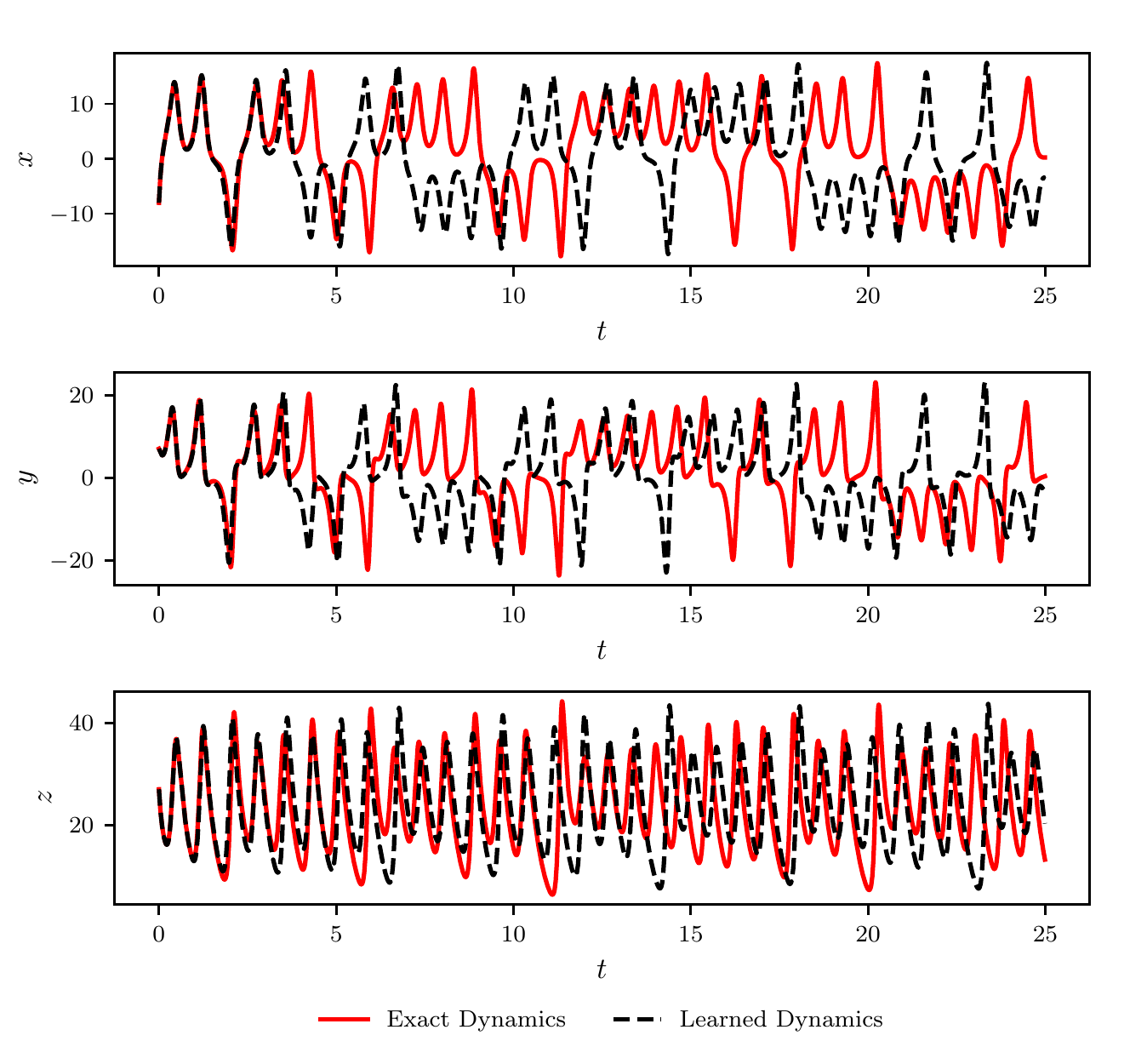}
\caption{{\em Lorenz System:} The exact trajectories of the Lorenz systems is compared to the corresponding trajectories of the learned dynamics. Solid blue lines represent the exact dynamics while the dashed black lines demonstrate the learned dynamics.}
\label{fig:Lorenz_Traj}
\end{figure}

To explore the identification of chaotic dynamics evolving on a finite dimensional attractor, we consider the nonlinear Lorenz system \cite{lorenz1963deterministic}
\begin{equation}
\begin{array}{l}
\dot{x} = 10 (y - x),\\
\dot{y} = x (28 - z) - y,\\
\dot{z} = x y - (8/3) z.
\end{array}
\end{equation}
We use $[x_0\ y_0\ z_0]^T = [-8\ 7\ 27]^T$ as initial condition and collect data from $t = 0$ to $t = 25$ with a time-step size of $\Delta t = 0.01$. The data are plotted in figures \ref{fig:Lorenz} and \ref{fig:Lorenz_Traj}. We employ a neural network with one hidden layer and $256$ neurons to represent the nonlinear dynamics. As for the multistep scheme \eqref{eq:multistep} we use Adams-Moulton with $M=1$ steps (i.e., the trapezoidal rule). Upon training the neural network, we solve the identified system using the same initial condition as the one above. As depicted in figure \ref{fig:Lorenz}, the learned system correctly captures the form of the attractor. \\

The Lorenz system has a positive Lyapunov exponent, and small differences between the exact and learned models grow exponentially, even though the attractor remains intact. This behavior is evident in figure \ref{fig:Lorenz_Traj}, as we compare the exact versus the predicted trajectories. Small discrepancies due to finite accuracy in the predicted dynamics lead to large errors in the forecasted time-series after $t>4$, despite the fact that the bi-stable structure of the attractor is well captured (see figure \ref{fig:Lorenz}).

\subsection{Fluid flow behind a cylinder}\label{sec:NavierStokes}

\begin{figure}[!t]
\includegraphics[width = 1.0\textwidth]{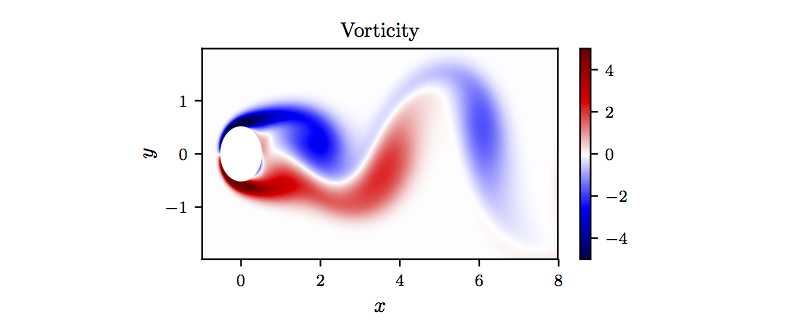}
\caption{{\em Flow past a cylinder:} A snapshot of the vorticity field of a solution to the Navier-Stokes equations for the fluid flow past a cylinder.}
\label{fig:Cylinder_vorticity}
\end{figure}

\begin{figure}[!t]
\includegraphics[width = 1.0\textwidth]{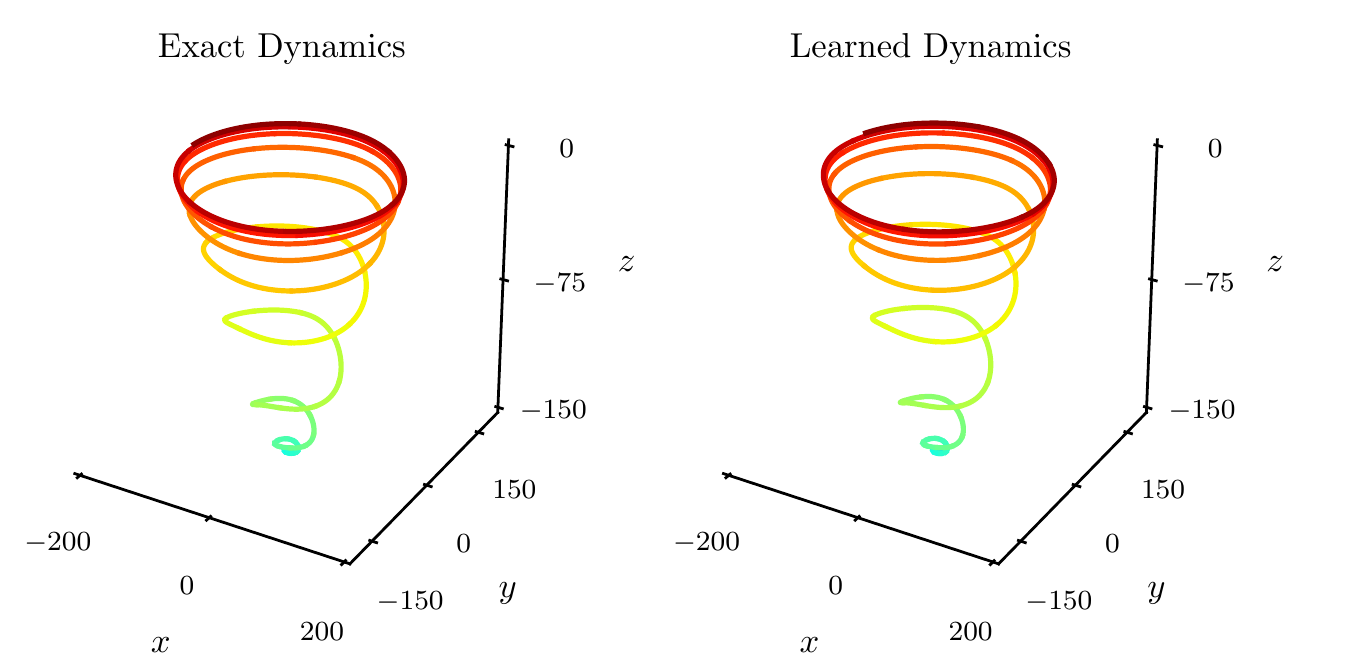}
\caption{{\em Flow past a cylinder:} The exact phase portrait of the cylinder wake trajectory in reduced coordinates (left panel) is compared to the corresponding phase portrait of the learned dynamics (right panel).}
\label{fig:Cylinder}
\end{figure}

In this example we collect data for the fluid flow past a cylinder (see figure \ref{fig:Cylinder_vorticity}) at Reynolds number 100 using direct numerical simulations of the two dimensional Navier-Stokes equations. In particular, following the problem setup presented in \cite{kutz2016dynamic} and \cite{Rudye1602614}, we simulate the Navier-Stokes equations describing the two-dimensional fluid flow past a circular cylinder at Reynolds number 100 using the Immersed Boundary Projection Method \cite{taira2007immersed, colonius2008fast}. This approach utilizes a multi-domain scheme with four nested domains, each successive grid being twice as large as the previous one. Length and time are non-dimensionalized so that the cylinder has unit diameter and the flow has unit velocity. Data is collected on the finest domain with dimensions $9 \times 4$ at a grid resolution of $449 \times 199$. The flow solver uses a 3rd-order Runge Kutta integration scheme with a time step of t = 0.02, which has been verified to yield well-resolved and converged flow fields. After simulations converge to steady periodic vortex shedding, flow snapshots are saved every $\Delta t = 0.02$. We then reduce the dimension of the system by proper orthogonal decomposition (POD) \cite{berkooz1993proper,brunton2016discovering}. The POD results in a hierarchy of orthonormal modes that, when truncated, capture most of the energy of the original system for the given rank truncation. The first two most energetic POD modes capture a significant portion of the energy; the steady-state vortex shedding is a limit cycle in these coordinates \cite{brunton2016discovering}. An additional mode, called the shift mode, is included to capture the transient dynamics connecting the unstable steady state with the mean of the limit cycle \cite{noack2003hierarchy}. The resulting POD coefficients are depicted in figure \ref{fig:Cylinder}. \\

We employ a neural network with one hidden layer and $256$ neurons to represent the nonlinear dynamics shown in figure \ref{fig:Cylinder}. As for the linear multistep scheme \eqref{eq:multistep} we use Adams-Moulton with $M=1$ steps (i.e., the trapezoidal rule). Upon training the neural network, we solve the identified system. As depicted in figure \ref{fig:Cylinder}, the learned system correctly captures the form of the dynamics and accurately reproduces the phase portrait, including both the transient regime as well as the limit cycle attained once the flow dynamics converge to the well known K\'{a}rman vortex street.

\subsection{Hopf bifurcation}\label{sec:Hopf}
\begin{figure}[!t]
\includegraphics[width = 1.0\textwidth]{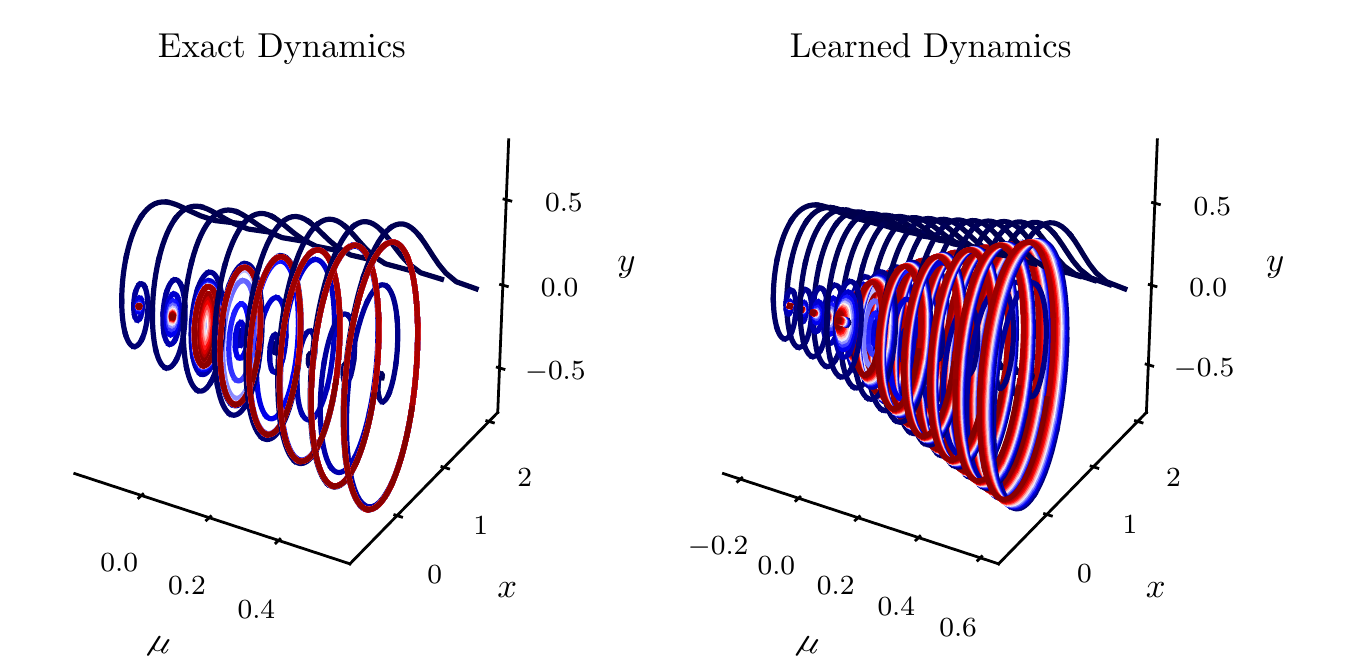}
\caption{{\em Hopf bifurcation:} Training data from the Hopf system for various initial conditions corresponding to different parameter values for $\mu$ (left panel) is compared to the corresponding phase portrait of the learned dynamics (right panel). It is worth highlighting that the algorithm is tested on initial conditions different from the ones used during training.}\label{fig:Hopf}
\end{figure}

Many real-world systems depend on parameters and, when the parameters are varied, they may go through bifurcations. To illustrate the ability of our method to identify parameterized dynamics, let us consider the Hopf normal form
\begin{equation}\label{eq:Hopf-normal}
\begin{array}{l}
\dot{x} = \mu x + y - x(x^2 + y^2),\\
\dot{y} = -x + \mu y - y(x^2 + y^2).
\end{array}
\end{equation}
Our algorithm can be readily extended to encompass parameterized systems. In particular, the system \eqref{eq:Hopf-normal} can be equivalently written as
\begin{equation}\label{eq:Hopf}
\begin{array}{l}
\dot{\mu} = 0,\\
\dot{x} = \mu x + y - x(x^2 + y^2),\\
\dot{y} = -x + \mu y - y(x^2 + y^2).
\end{array}
\end{equation}
We collect data from the Hopf system \eqref{eq:Hopf} for various initial conditions corresponding to different parameter values for $\mu$. The data is depicted in figure \ref{fig:Hopf}. The identified parameterized dynamics is shown in figure \ref{fig:Hopf} for a set of parameter values different from the ones used during model training. The learned system correctly captures the transition from the fixed point for $\mu < 0$ to the limit cycle for $\mu>0$.

\subsection{Glycolytic oscillator}\label{sec:glycolytic}
\begin{figure}[!t]
\includegraphics[width = 1.0\textwidth]{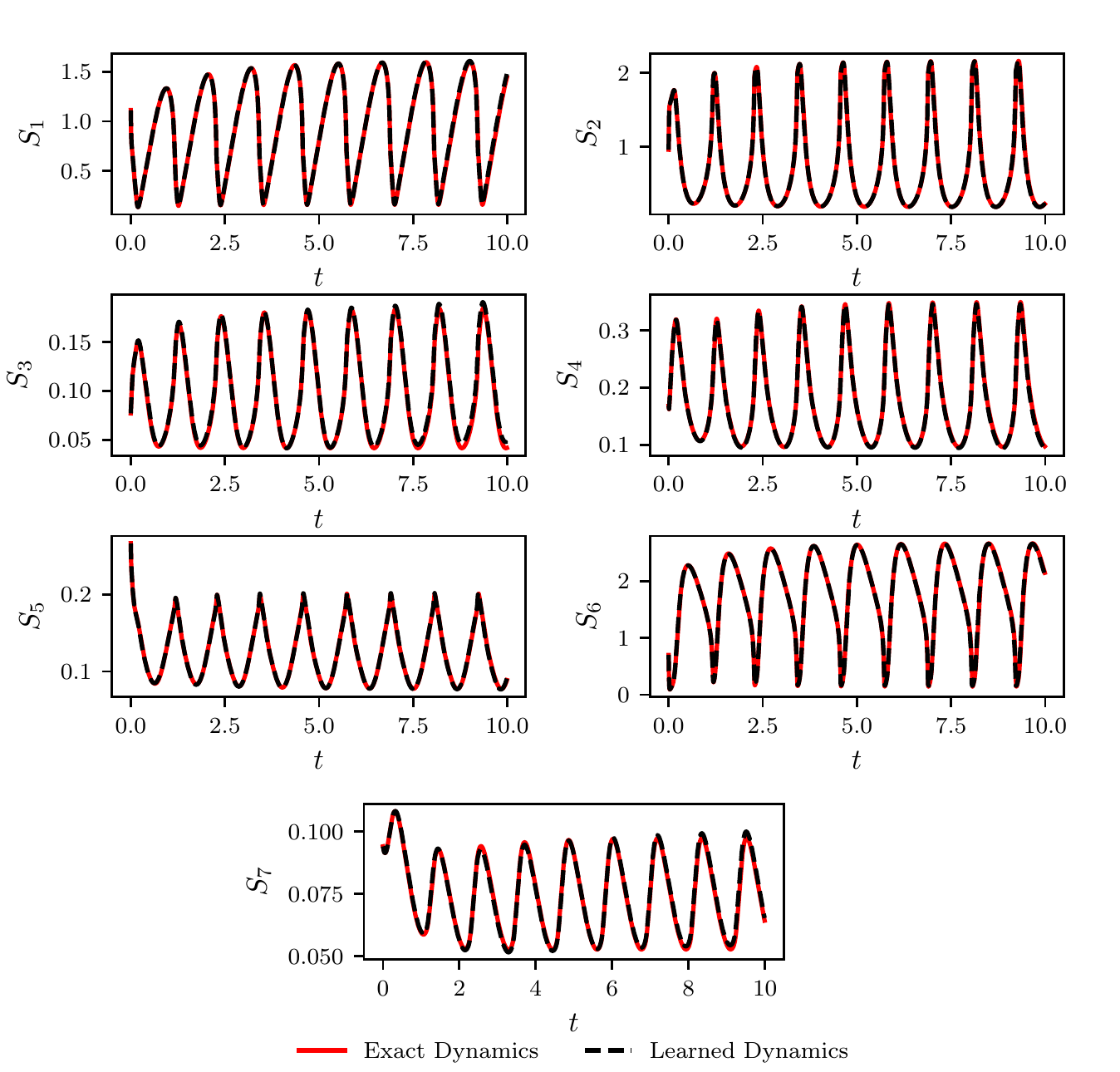}
\caption{{\em Glycolytic oscillator:} Exact versus learned dynamics for random initial conditions chosen from the ranges provided in table 2 of \cite{daniels2015efficient}.}
\label{fig:Glycolytic}
\end{figure}

As an example of complicated nonlinear dynamics typical of biological systems, we simulate the glycolytic oscillator model presented in \cite{daniels2015efficient} and \cite{brunton2016discovering}. The model consists of ordinary differential equations for the concentrations of 7 biochemical species; i.e.,
\begin{eqnarray}\label{eq:glycolytic}
&&\frac{dS_1}{dt} = J_0 - \frac{k_1 S_1 S_6}{1 + (S_6/K_1)^q},\nonumber\\
&&\frac{dS_2}{dt} = 2\frac{k_1 S_1 S_6}{1 + (S_6/K_1)^q} - k_2 S_2 (N - S_5) - k_6 S_2 S_5,\nonumber\\
&&\frac{dS_3}{dt} = k_2 S_2 (N - S_5) - k_3 S_3 (N - S_6),\nonumber\\
&&\frac{dS_4}{dt} = k_3 S_3 (A - S_6) - k_4 S_4 S_5 - \kappa (S_4 - S_7),\\
&&\frac{dS_5}{dt} = k_2 S_2 (N - S_5) - k_4 S_4 S_5 - k_6 S_2 S_5,\nonumber\\
&&\frac{dS_6}{dt} = -2\frac{k_1 S_1 S_6}{1 + (S_6/K_1)^q} + 2 k_3 S_3 (A - S_6) - k_5 S_6,\nonumber\\
&&\frac{dS_7}{dt} = \psi \kappa (S_4 - S_7) - k S_7.\nonumber
\end{eqnarray}
The parameters of the model are chosen according to table 1 of \cite{daniels2015efficient}. As shown in figure \ref{fig:Glycolytic}, data from a simulation of equation \eqref{eq:glycolytic} are collected from $t = 0$ to $t = 10$ with a time-step size of $\Delta t = 0.01$. We employ a neural network with one hidden layer and $256$ neurons to represent the nonlinear dynamics. As for the multi-step scheme \eqref{eq:multistep} we use Adams-Moulton with $M=1$ steps (i.e., the trapezoidal rule). Upon training the neural network, we solve the identified system using the same initial condition as the ones used for the exact system. As depicted in figure \ref{fig:Glycolytic}, the learned system correctly captures the form of the dynamics.

\section{Summary and Discussion}\label{sec:Conclusion} 

We have presented a machine learning approach for extracting nonlinear dynamical systems from time-series data. The proposed algorithm leverages the structure of well studied multi-step time-stepping schemes such as Adams-Bashforth, Adams Moulton, and BDF families, to construct efficient algorithms for learning dynamical systems using deep neural networks. A key property of the proposed approach is the use of multiple steps which enables us to incorporate memory effects in learning the temporal dynamics and tackle problems with a nonlinear and non-Markovian dynamical structure. Specifically, the use of $M$ steps allows us to decouple the regression complexity due to several temporal lags, ultimately leading to a simpler $D$-dimensional regression problem, as opposed to an $(M\times D)$-dimensional problem in the case of a brute force NARMAX or recurrent neural network approaches. Although state-of-the-art results are presented for a diverse collection of benchmark problems, there exist a series of open questions mandating further investigation. How could one handle a variable temporal gap $\Delta{t}$, i.e., irregularly sampled data in time? How would common techniques such as batch normalization, drop out, and $\mathcal{L}_1$/$\mathcal{L}_2$ regularization enhance the robustness of the proposed algorithm and mitigate the effects of over-fitting? How could one incorporate partial knowledge of the dynamical system in cases where certain interaction terms are already known? In terms of future work, interesting directions include the application of convolutional architectures \cite{goodfellow2016deep} for mitigating the complexity associated with very high-dimensional inputs, as well as studying possible connections with recent studies linking deep neural networks with numerical methods and dynamical systems \cite{chang2017multi,lu2017beyond}.  

\section*{Acknowledgements}
This work received support by the DARPA EQUiPS grant N66001-15-2-4055 and the AFOSR grant FA9550-17-1-0013. All data and codes used in this manuscript are publicly available on GitHub at \url{https://github.com/maziarraissi/MultistepNNs}.

% \appendix
% \input{appendix.tex}


%% References


\bibliographystyle{model1-num-names}
\bibliography{sample.bib}






% \end{thebibliography}

\end{document}